\documentclass[10pt]{article}
\pdfoutput=1 
\usepackage{a4}
\usepackage{amsmath}
\usepackage{amssymb}
\usepackage{amsthm}
\usepackage{lipsum} 
\usepackage{amsfonts}
\usepackage{graphicx}
\usepackage{epstopdf}
\usepackage{epsfig}
\usepackage{pgfgantt}
\usepackage[numbers, sort&compress]{natbib}
\usepackage{cite}
\usepackage{enumerate}
\usepackage{parskip}
\usepackage{cleveref}
\usepackage[affil-it]{authblk}

\newcommand{\EQ}[2]{\begin{equation}{{#2}\label{#1}} \end{equation}}
                        \newcounter{theorem}

\renewcommand{\thetheorem}{\arabic{section}.\arabic{theorem}}

\newenvironment{thm}[2]{\begin{sloppypar}\refstepcounter{theorem}%
                        {\bf #1 \thetheorem.}\label{#2}\em{}}%
                        {\end{sloppypar}}
                        \newcommand{\theo}[3]{\begin{thm}{#1}{#2}
#3\end{thm}}

                        \newcommand{\N}{{\rm I}\!{\rm N}}
             
                         \newcommand{\Z}{\makebox[0.04cm][l]{\sf Z}\!{\sf Z}}

                        \newcommand{\R}{{\rm I}\!{\rm R}}
                        \newcommand{\C}{\makebox[0.15cm][l]{\rm
C}\!{\sf I}\hspace{0.13cm}}

\def\endofproof{\vbox{\hrule\hbox{\vrule\hbox to 6truept
                                      {\vbox to 6truept
{\vfil}\hfil}\vrule}\hrule}}%

\title{Generalised Wendland functions for the sphere}

\author[1]{Simon Hubbert}
\author[2]{ Janin J\"ager\thanks{corrsponding author: janin.jaeger@math.uni.giessen.de}}
\affil[1]{School of Economics, Mathematics and Statistics, Birkbeck, University of London, Malet Street, London, WC1E 7HX,United Kingdom}
\affil[2]{Mathematisches Institut,
Justus-Liebig University, 35392 Giessen, Germany.}

\begin{document}
\def\R{{\mathbb R}}
\def\Z{{\mathbb Z}}
\def\N{{\mathbb N}}
\def\S{{\mathbb{S}^2}}
\def\wh{\widehat}
\def\C{{\mathbb C}}
\def\bigO{{\mathcal{O}}}

\maketitle

\begin{abstract}
In this paper we compute the spherical Fourier expansions coefficients for
 the restriction of the generalised Wendland functions from $d-$dimensional Euclidean space to the $(d-1)-$dimensional unit sphere. The development required to derive these coefficients relies heavily upon known asymptotic results for hypergeometric functions and the final result shows that they can be expressed in closed form as a multiple of a certain $_{3}F_{2}$ hypergeometric function. Using the closed form expressions we are able to provide the precise asymptotic rates of decay for the spherical Fourier coefficients which we observe have a close connection to the asymptotic decay rate of the corresponding Euclidean Fourier transform.
\end{abstract}

\textit{keywords:} positive definite kernels, spherical basis functions, compact support. 

\section{Introduction}

Positive definite functions are frequently used in
 scattered data fitting algorithms both in Euclidean space and on
 spheres: see \citep{WendBook}. The special case of the $2$-sphere is of importance in geostatistics where positive definite functions are used as covariance functions of random fields on the surface of the earth,  \citep{Lang2015}. The aim of this paper is to
 investigate the generalised Wendland functions, a parameterised family of compactly supported basis functions  which, for a certain parameter range, are strictly  positive definite on $\R^d$ . In the opening section we will review the rudimentary material relating to positive definite kernels defined on  Euclidean space and on the unit sphere. We will pay particular attention to the restriction of Euclidean positive definite functions to the sphere and flag a crucial identity that connects the Fourier transform of the Euclidean function to the spherical Fourier coefficients of its restriction to the sphere. In section 3 we introduce the family of  generalised Wendland functions, setting out their known properties including an expression for their Euclidean Fourier transforms. Then, in section 4,  we make use of the aforementioned connection identity to derive an expression involving a certain hypergeometric function for the spherical Fourier coefficients of the generalised Wendland functions restricted to the sphere. In the final section we derive the precise asymptotic behaviour of the spherical Fourier coefficients which, as we will show, is very closely linked to the analogous asymptotic decay of the Euclidean Fourier transform.
 
The class of generalised Wendland functions, as their name suggests, contains the original Wendland functions, which are popular in applications due to their simple polynomial form.  Many researchers have employed the original Wendland functions on the sphere and, for these functions, the asymptotic behaviour of the spherical Fourier coefficients has been addressed. However, as far as the authors are aware, a precise formula for the coefficients has not previously been made available. A closed form of the coefficients is necessary for the use of recently proposed numerical methods such as the stable computation via  Hilbert-Schmidt SVD (\citep{Fass}, Chapter 13) and also for the spectral simulation of Gaussian random fields as described in  \citep{Lantuejoul2019}.

Of the work that is related to ours we draw the reader's attention to \citep{MNPW2010} (Proposition 3.1)  where a precise asymptotic form for the spherical Fourier coefficients of the original Wendland functions is derived, however the constant multiplying the decay factor is not explicitly given.  In addition,  le Gia et al. \citep{LGSW} (Section 6) have consider  scaled versions of positive definite functions and show that if their Fourier transforms decay at a polynomial rate then, when restricted to the sphere,  their corresponding Fourier coefficients  decay at the expected analogous rate; this work applies to the class of functions under consideration in this paper however only upper and lower bounds on the rate of decay can be inferred. In our work, by deriving a precise formula for the Fourier coefficients of a much broader family of functions we are also able to determine their precise asymptotic decay rates.
 
\section{Radial and zonal kernels}

\theo{Definition}{pdrad}{A  kernel $\Phi:\Omega \times \Omega \to \R$ is said to
be strictly positive definite on a domain 
$\Omega,$ if, for any $n\ge 2$ distinct locations $x_{1},\ldots,x_{n}
\in \Omega,$ the  $n \times n$  matrix
\EQ{distmat}{\Bigl(\Phi(x_{j},x_{k})\Bigr)_{j,k = 1}^{n}, }
is symmetric and  positive definite.
}


For the class of radial kernels taking the form $\Phi(x,y)=\phi(\|x-y\|)$ we have the following  characterisation theorem
(see \citep{WendBook} Theorem 6.18).

\theo{Theorem}{char}{A radial kernel $\Phi(x,y)=\phi(\|x-y\|)$, with  $\phi: [0,\infty)\to \R$
such that $r \mapsto r^{d-1}\phi(r) \in L_{1}[0,\infty)$ is strictly
 positive definite  on $\R^{d}$ if and only
if the $d-$dimensional Fourier transform \begin{equation}
 \widehat{\phi}(z)=
z^{1-\frac{d}{2}} \int_0^{\infty} \phi(y) \, y^{\frac{d}{2}} \,
J_{\frac{d}{2}-1}(yz) \, dy, \label{eqnRadialFourierTrans}
\end{equation}
(where $J_{\nu}(\cdot)$ denotes the Bessel function of the first
kind with order $\nu$) is non-negative and not identically equal to
zero.}


If we assume that  $\widehat{\phi}(z)>0$ then we can appeal to the theory of radial basis functions (see \citep{WendBook})  to deduce the following result. 

\theo{Theorem}{RBFnativeSpaces}{ Let $d\ge 1$ denote a fixed
spatial dimension and $\Phi$ be a strictly positive definite radial kernel on $\R^{d}$ with $\widehat{\phi}(z)>0$ for all $z \ge 0.$ 
Define
\EQ{native}{
N_{\phi}:=\Big\{f \in L_{2}(\R^{d}): \|f\|_{\phi}^{2}=\frac{1}{(2\pi)^{\frac{d}{2}}}\int_{\R^{d}}\frac{|\widehat{f}(\boldsymbol{\omega})|^{2}}
{\widehat{\phi}(\|\boldsymbol{\omega}\|)}d\boldsymbol{\omega}\Big\}
}
where $\|\cdot\|_{\phi}$ is a norm induced by the inner product
\EQ{inner}
{(f,g)_{\phi}:=\frac{1}{(2\pi)^{\frac{d}{2}}}\int_{\R^{d}}\frac{\widehat{f}(\boldsymbol{\omega})\widehat{g}(\boldsymbol{\omega})}
{\widehat{\phi}(\|\boldsymbol{\omega}\|)}d\boldsymbol{\omega}.
}
Then $N_{\phi}$ is a real Hilbert space with inner product $(\cdot,\cdot)_{\Phi}$ and reproducing kernel  $\Phi$. }


\theo{Remark}{RBFsob}{ 
We note that if there exists positive constants $c_1<c_2$ such that
\EQ{decay}{
\frac{c_{1}}{(1+z^{2})^\lambda}\le \widehat{\phi}(z)\le  \frac{c_{2}}{(1+z^{2})^\lambda}, \quad z\in \R,
}
where $\lambda >\frac{d}{2},$ then
the space $N_{\phi}$ (\ref{native}) is norm equivalent to the Sobolev space
$$
H^{\lambda}(\R^{d}):=\Bigl\{f \in L_{2}(\R^{d})\cap C(\R^{d}): \int_{\R^{d}}|\widehat{f}(\boldsymbol{\omega})|^{2}(1+\|\boldsymbol{\omega}\|^{2})^{\lambda}< \infty \Big\}.
$$
}

We now consider the spherical case. We let $S^{d-1}:=\{ {x} \in \R^{d}:\|x\|=1\},$ denote the $(d-1)-$dimensional unit sphere, then for any points $\boldsymbol\xi,\boldsymbol\eta \in S^{d-1}$ we write $\boldsymbol\xi^{T} \boldsymbol\eta = \cos(\theta)$ to denote their dot-product  where $\theta$ is the angular distance between the two points. 


The zonal kernel $\Psi(\boldsymbol\xi,\boldsymbol\eta)=\psi(\boldsymbol\xi^{T}\boldsymbol\eta)$ induced by any continuous $\psi:[-1,1]\to \R,$
 possesses  a Fourier-type expansion in spherical harmonics
\EQ{expansion}{
\Psi(\boldsymbol\xi,\boldsymbol\eta) =\psi(\boldsymbol{\xi}^{T}\boldsymbol{\eta})=\sum\limits_{m=0}^{\infty}\sum\limits_{n=1}^{N_{m,d}}\widehat{\psi}_{m} {\cal{Y}}_{m,n}(\boldsymbol\xi){\cal{Y}}_{m,n}(\boldsymbol\eta).
}
where $\{{\cal{Y}}_{m,n}:n=1,\ldots, N_{m,d}\}$ is an orthonormal basis for the space of spherical harmonics of degree $m$  and the collection $\{{\cal{Y}}_{m,n}:n=1,\ldots, N_{m,d}, m\ge 0\}$
forms an orthonormal basis for $L_{2}(S^{d-1}).$

Using Schoenberg's \citep{Schoen} pioneering work  it can be shown that if the expansion coefficients $\widehat{\psi}_{m} $ (commonly referred to as the Schoenberg coefficients) are strictly positive for  $m\ge0$ then $\Psi$ is a positive definite kernel on $S^{d-1}.$ This simple
condition, as we shall see, is sufficient for our purposes but the reader may consult \citep{CheMenSun03} for a careful investigation of the necessary \textsl{and} sufficient conditions.

Spherical harmonics provide a Fourier analysis for the sphere. In particular, every $f \in L_{2}(S^{d-1})$ has an associated spherical Fourier expansion
\[
f = \sum_{m=0}^{\infty}\sum_{n=1}^{N_{n,d}}\widehat{f}_{n,m}{\cal{Y}}_{m,n}\,\,\,\,{\rm{where}}\,\,\,\,
\widehat{f}_{n,m}=(f,{\cal{Y}}_{m,n})_{L_{2}(S^{d-1})}.
\]
The following theorem is  the spherical analogue of Theorem \ref{RBFnativeSpaces}.

\theo{Theorem}{SBFnativeSpaces}{ Let $d\ge 2$ denote a fixed
spatial dimension and $\Psi$  a strictly positive definite zonal kernel  on $S^{d-1}$ for which the Fourier expansion coefficients $\widehat{\psi}_{m}$ are strictly positive for all $m \ge 0.$ Define

\EQ{natspsph}{N_{\psi}=\Bigl\{f\in L_{2}(S^{d-1}):\|f\|_{\psi}^{2}=\sum_{m=0}^{\infty}\sum_{n=1}^{N_{n,d}}\frac{|\widehat{f}_{n,m}|^{2}}{\widehat{\psi}_{m}} <\infty \Bigr\}}
where $\|\cdot\|_{\psi}$ is a norm induced by the inner-product
\[
(f,g)_{\psi}:=\sum_{m=0}^{\infty}\sum_{n=1}^{N_{n,d}}\frac{\widehat{f}_{n,m}\widehat{g}_{n,m}}{\widehat{\psi}_{m}}.
\]
Then $N_{\psi}$ is a real Hilbert space with inner product $(\cdot,\cdot)_{\Phi}$ and reproducing kernel  $\Psi(\boldsymbol{\xi},\boldsymbol{\eta}).$ }



\theo{Remark}{SBFsob}{ 
We note that if there exists positive constants $c_1<c_2$ such that
\EQ{decaysp}{
\frac{c_{1}}{(1+m^{2})^\lambda}\le \widehat{\psi}_{m}\le  \frac{c_{2}}{(1+m^{2})^\lambda}, \quad m\ge  0,
}
where $\lambda >\frac{d-1}{2},$ then
the space $N_{\psi}$ defined in the previous theorem is norm equivalent to the Sobolev space
$$
H^{\lambda}(S^{d-1}):=\Bigl\{f \in L_{2}(S^{d-1})\cap C(S^{d-1}):\sum_{m=0}^{\infty}\sum_{n=1}^{N_{n,d}}(1+m^2)^{\lambda}|\widehat{f}_{n,m}|^{2} <\infty \Bigr\}.
$$
}

In practice one can select $\phi$ to be a positive definite function  on $\R^{d}$ and use the relation $\|\boldsymbol\xi - \boldsymbol\eta\| = \sqrt{2-2\boldsymbol\xi^{T}\boldsymbol\eta},$ for $\boldsymbol\xi,\boldsymbol\eta \in S^{d-1},$
to define its restriction to $S^{d-1}$
as
\EQ{restrict}{
\Psi(\boldsymbol\xi,\boldsymbol\eta)=\psi(\boldsymbol\xi^{T}\boldsymbol\eta)=\phi \left(\sqrt{2-2\boldsymbol\xi^{T}\boldsymbol\eta}\right).
}
In this regard we have the following formula (\citep{NW02} Theorem 4.1) which links the spherical Fourier coefficients of $\Psi$ to the radial Fourier transform  $\widehat{\phi}(z),$
\EQ{linkup}
{\widehat{\psi}_{m}:=(2\pi)^{\frac{d}{2}}\int_{0}^{\infty}z J_{m+\frac{d-2}{2}}^{2}(z)\widehat{\phi}(z)dz.}


\section{The generalised Wendland functions}

We will investigate a family of parameterised basis
functions that is generated by a  truncated power function. Specifically, we choose a support parameter $\epsilon >0,$ and define
\[
\phi_{\mu,0}^{(\epsilon)}(r):=(1-\epsilon r)_{+}^{\mu}=\begin{cases} (1-\epsilon r)^{\mu} \,\, & \,\,
\textrm{for} \,\,\,\, 0\le r\le \frac{1}{\epsilon}; \\
0 \,\, & \,\,
\textrm{for} \,\,\,\, r\ge \frac{1}{\epsilon},
\end{cases}
\]
and consider

\EQ{family}{ \phi_{\mu,\alpha}^{(\epsilon)}(r) :=
\frac{1}{2^{\alpha-1}\Gamma(\alpha)}\int_{\epsilon r}^{1} \phi_{\mu,0}(t) \, t \,
\left(t^2-(\epsilon r)^2 \right)^{\alpha-1} \mathrm{d}t \quad \mbox{for} \quad
r \in \left[0,\frac{1}{\epsilon}\right], } where $\mu > -1,$ $\alpha>0$ and $\Gamma(\cdot)$
denotes the Gamma function.

This class was also referred to as Zastavnyi functions and discussed in  \citep{Zastavnyi1}. It is a subclass of the Buhmann functions originally introduced in \citep{Buhmann}. Applications of this function class as covariance functions of Gaussian random fields are discussed in \citep{Bevilacqua}. 

We remark that with a judicious  selection  of the $\mu$ and $\alpha$ parameters we can recover, from formula (\ref{family}), both  the original, and most well-known, Wendland functions and also the so-called missing Wendland functions. It is for this reason that this family is often referred to as the generalised Wendland functions.  In more detail, the original Wendland functions are recovered when the space dimension is odd, $\alpha$ is a positive integer, i.e, $\alpha:=k \in \Z_{+},$ and $\mu:=\ell= \frac{d+1}{2}+k,$ i.e., the smallest integer 
that still allows positive definiteness. In this case one can show that 
\[
\phi_{\ell,k}^{(\epsilon)}(r)= p_{k}(\epsilon r)(1-\epsilon r)_{+}^{\ell+k},
\]
where $p_{k}$ is a polynomial of degree $k.$ The missing Wendland functions are recovered when the space dimension is even,  $\alpha$ is a positive half-integer, i.e., $\alpha = k-\frac{1}{2}$ (where $k \in \Z_{+}),$  and $\mu:=\ell= \frac{d}{2}+k,$ i.e., once again the smallest integer 
that still allows positive definiteness. The missing Wendland functions have 
two polynomial components, one with a logarithmic multiplier $L(r):=\log\left(\frac{r}{1+\sqrt{1-r^{2}}}\right) $ and one with a square root multiplier $S(r):=\sqrt{1-r^{2}},$ see \citep{SchabMiss}.

The $d$-dimensional Fourier transform of $\phi_{\mu,\alpha}^{(\epsilon)}$ was computed in  \citep{CH} and it can be expressed succinctly as  
\EQ{FTeuc}{
\begin{aligned}
\widehat{\phi_{\mu,\alpha}^{(\epsilon)}}(z)=& \frac{C_{\lambda,\mu}}{\sqrt{2\pi}\epsilon^{d}}   {}_{1}F_2\left( \lambda; \lambda+ \frac{\mu}{2},\lambda+ \frac{\mu+1}{2}; -\left(\frac{z}{2\epsilon}\right)^2\right),
\end{aligned}}
where
\EQ{ClamMu}{
 \lambda := \frac{d+1}{2}+\alpha\quad {\rm{and}}\quad C_{\lambda,\mu}=\frac{2^{\lambda}\Gamma(\lambda)\Gamma(\mu+1)}{\Gamma(2\lambda +\mu)}}
 and
 where ${}_{1}F_2(a;b,c;z)$ denotes the
hypergeometric
    function (see \citep{Abramowitz1964}, (15.1.1)). Hypergeometric functions will feature heavily in this work and so we  briefly remind the reader that a general
hypergeometric
    function is
   defined by
    \EQ{hyper}{
_{p}F_{q}(a_1,\ldots,a_p;b_1,\ldots,b_q;z):=\sum_{j=0}^{\infty}\frac{\Pi_{i=1}^p
(a_i)_j }{\Pi_{i=1}^q (b_i)_j}\frac{z^{j}}{j!}, }
where \EQ{poch}{ (c)_n := c(c+1)\cdots(c+n-1) =
\frac{\Gamma(c+n)}{\Gamma(c)}, \hspace{0.1in} n \geq 1} denotes the Pochhammer symbol, with $(c)_0 = 1$.

It is known (see \citep{CH}) that $\widehat{\phi_{\mu,\alpha}^{(\epsilon)}}(z)>0$  if and only if its parameters satisfy $
\mu \ge \lambda.$ Thus, with such a choice,  $ \phi_{\mu,\alpha}^{(\epsilon)}(r)$ induces a strictly positive definite and compactly supported radial kernel on $\R^{d}.$  The following formula, taken from \citep{Miller1998} (equation 2.2), provides the asymptotic behaviour 
for  ${}_1F_2$ hypergeometric functions for large argument $t$
\[
\begin{aligned}
{}_1F_2\left( \alpha; \beta,1+\nu
 ;\ -t^2\right)&=\frac{\Gamma\left(\beta\right)\Gamma\left(1+\nu\right)}{\Gamma\left(\beta-\alpha\right)\Gamma\left(1+\nu-\alpha\right)}\frac{1}{t^{2\alpha}}\\
 &+
 \frac{\Gamma\left(\beta\right)\Gamma\left(1+\nu\right)}{\sqrt{\pi}\Gamma\left(\alpha\right)}\frac{\cos\Bigl[2t-\frac{\pi}{2}\left(\beta+\nu+\frac{1}{2}-\alpha\right)\Bigr]}{t^{\beta+\nu+\frac{1}{2}-\alpha}}.
 \end{aligned}
 \]
 
 Setting $\alpha:=\lambda,$ $\beta:=\lambda+\frac{\mu}{2},$ $\nu:=\lambda+\frac{\mu-1}{2}$ and $t:=\frac{z}{2\epsilon}$ we deduce that
 \EQ{ass1F2}{
 \begin{aligned}
& {}_1F_2\left( \lambda; \lambda+\frac{\mu}{2},\lambda+\frac{\mu+1}{2}
 ;\ -\frac{z^2}{4\epsilon^2}\right)= \Gamma\left(\lambda+\frac{\mu}{2}\right)\Gamma\left(\lambda+\frac{\mu+1}{2}\right) \\
 &
\times\left[\frac{\left(\frac{2\epsilon}{z}\right)^{2\lambda}}{\Gamma\left(\frac{\mu}{2}\right)\Gamma\left(\frac{\mu+1}{2}\right)}
 +\frac{\cos\Bigl[\frac{z}{\epsilon}
 -\frac{\pi}{2}\left(\lambda+\mu\right)\Bigr]\left(\frac{2\epsilon}{z}\right)^{\lambda+\mu}}{\sqrt{\pi}\Gamma\left(\lambda\right)}
 \right].
 \end{aligned}
 } 
 
Since $\mu\ge \lambda>\frac{1}{2}$, it is the first term that determines the asymptotic decay and we can conclude, from \eqref{FTeuc},  that
\EQ{assdecayR}{
\begin{aligned}
\widehat{\phi_{\mu,\alpha}^{(\epsilon)}}(z)&\sim 
 \frac{C_{\lambda,\mu}}{\sqrt{2\pi}\epsilon^{d}} \frac{\Gamma\left(\lambda+\frac{\mu}{2}\right)\Gamma\left(\lambda+\frac{\mu+1}{2}\right)}{\Gamma\left(\frac{\mu}{2}\right)\Gamma\left(\frac{\mu+1}{2}\right)}\left(\frac{2\epsilon}{z}\right)^{2\lambda}=\frac{2^{\lambda}\Gamma(\lambda)\mu}{\sqrt{2\pi}}\frac{\epsilon^{2\alpha+1}}{z^{2\lambda}},\\
\end{aligned}
}
where the equality follows from applying the Gamma function duplication formula  (\citep{Abramowitz1964}, (6.1.18))
\EQ{gammadup}{
\frac{\Gamma(2z)}{\Gamma\left(z+\frac{1}{2}\right)}=\frac{2^{2z-1}}{\sqrt{\pi}}\Gamma(z),
}
and recalling the definition of $\lambda$ and the constant $C_{\lambda,\mu}$
 (\ref{ClamMu}). This  asymptotic behaviour implies that there 
 exist positive constants $c_{1}<c_{2}$ such that
\EQ{decayGW}{
\frac{c_{1}\epsilon^{2\alpha+1}}{(1+z^{2})^\lambda}\le \widehat{\phi_{\mu,\alpha}^{(\epsilon)}}(z)\le  \frac{c_{2}\epsilon^{2\alpha+1}}{(1+z^{2})^\lambda}, \quad z\in \R.
}

Appealing to Theorem \ref{RBFnativeSpaces} and Remark \ref{RBFsob} we can  deduce that when $\mu \ge \lambda,$   $\Phi(\textbf{x},\textbf{y})  =
\phi_{\mu,\alpha}^{(\epsilon)}(\|\textbf{x}-\textbf{y}\|)$  defines a reproducing kernel of the Sobolev
space $H^{\lambda}(\R^{d}).$ 

\section{Generalized Wendland functions for the sphere}
In this section we will consider the restriction of $\phi_{\mu,\alpha}^{(\epsilon)}(r)$ to the unit sphere. Using (\ref{restrict})  and (\ref{expansion}) we can write
 \EQ{reskernel}{
 \phi_{\mu,\alpha}^{(\epsilon)}\left(\sqrt{2-2\boldsymbol\xi^{T}\boldsymbol\eta}\right)=\Psi_{\mu,\alpha}^{(\epsilon)}(\boldsymbol\xi,\boldsymbol\eta)=
 \sum\limits_{m=0}^{\infty}\sum\limits_{n=1}^{N_{m,d}}\widehat{\psi_{\mu,\alpha}^{(\epsilon)}}(m){\cal{Y}}_{m,n}(\boldsymbol\xi){\cal{Y}}_{m,n}(\boldsymbol\eta).
 }
 
 We note that the support condition of the restriction can be recast in terms of angular distance $\theta$ (between $ \boldsymbol\xi,\boldsymbol\eta \in S^{d-1}$) as 
 \[
0\le \sqrt{2-2\cos(\theta)}\le\frac{1}{\epsilon}\implies 1-\frac{1}{2\epsilon^2}\le\cos(\theta) \le 1.
\]
Thus, we need only consider the range $\epsilon \ge \frac{1}{2}$ for the support parameter; the case $\epsilon =\frac{1}{2}$ ensures that the restricted function is globally supported on the entire sphere, whereas $\epsilon>\frac{1}{2}$ ensures that it is supported on a spherical cap of radius $\theta = \cos^{-1}\left(1-\frac{1}{2\epsilon^2}\right).$


\theo{Theorem}{Coeff}{
The spherical Fourier coefficients ($d$-Schoenberg coefficients) of the generalised Wendland function $\Psi_{\mu,\alpha}^{(\epsilon)}(\boldsymbol\xi,\boldsymbol\eta)$ from \eqref{reskernel} are given by 
\EQ{swc}{
\begin{aligned}
 &\widehat{\psi_{\mu,\alpha}^{(\epsilon)}}(m)=(2\pi)^{\frac{d-1}{2}}\frac{C_{\lambda-\frac{1}{2},\mu}}{\sqrt{2\pi}\epsilon^{d-1}}\\
& \times
{}_3F_2\left(- \left( m+\frac{d-3}{2}\right),m+\frac{d-1}{2},\lambda-\frac{1}{2};\lambda+\frac{\mu-1}{2},\lambda+\frac{\mu}{2};
\frac{1}{4\epsilon^2}\right),
\end{aligned}
} where, in analogy with (\ref{ClamMu}), we have that
\EQ{ClamsphMu}{
 \lambda := \frac{d+1}{2}+\alpha\quad {\rm{and}}\quad C_{\lambda-\frac{1}{2},\mu}=\frac{2^{\lambda-\frac{1}{2}}\Gamma\left(\lambda-\frac{1}{2}\right)\Gamma(\mu+1)}{\Gamma(2\lambda +\mu-1)}}
}

We observe that the expression (\ref{swc}) shows a close connection to that of the Euclidean Fourier transform (\ref{FTeuc})  of the generalised Wendland functions.

To establish the result we consider, for $\gamma>0$,  the following finite integral
 \EQ{FinInt}{
 I_{\gamma}(r) = (2\pi)^{\frac{d}{2}}\int_{0}^{r}z J_{\gamma}^{2}(z)\widehat{\phi_{\mu,\alpha}^{(\epsilon)}}(z)dz.}
We note that as a specific instance of  the connection formula (\ref{linkup}) the associated spherical Fourier coefficients are given by
\EQ{GWcoefs}{
 \widehat{\psi_{\mu,\alpha}^{(\epsilon)}}(m)= \lim_{r\to \infty}I_{m+\frac{d-2}{2}}(r).
 }
 
In order to prove the above theorem we determine a closed form representation of this finite integral.

\theo{Lemma}{FinitInt}{The integral defined in \eqref{FinInt} satisfies:
\EQ{first}{
\begin{aligned}
I_{\gamma}(r)=&(2\pi)^{\frac{d}{2}}
\frac{C_{\lambda,\mu}\left(2^{\gamma}\Gamma\left(\gamma+1\right)\right)^{-2} }{\sqrt{2\pi}\epsilon^{d}}\frac{r^{2(\gamma+1)}}{2(\gamma+1)}\sum_{\ell =0}^{\infty} \frac{(\gamma+1)_{\ell}(\lambda)_{\ell}\left(-\frac{r^2}{4\epsilon^2}\right)^{\ell}}{(\gamma+2)_{\ell}\left(\lambda+\frac{\mu}{2}\right)_{\ell}\left(\lambda+\frac{\mu+1}{2}\right)_{\ell}\ell!}
\\
&\times{}_{2}F_3\left(\gamma+1+\ell,\gamma+\frac{1}{2};\gamma+2+\ell,\gamma+1,2\gamma+1;-r^{2}\right).
\end{aligned}
}}

\begin{proof}

From \citep{Luke} (Section 6.2, (41)) we have that
\[
\left( J_{\gamma}(z)\right)^2=\left(2^{\gamma}\Gamma(\gamma+1)\right)^{-2}z^{2\gamma} {}_1F_{2}\left(\gamma+\frac{1}{2}; \gamma+1, 2\gamma+1; -z^2\right),
\]

and using  this and (\ref{FTeuc}) we can write the integral  (\ref{FinInt}) as

\EQ{gammaintegral}{
I_{\gamma}(r)=(2\pi)^{\frac{d}{2}}\frac{C_{\lambda,\mu}\left(2^{\gamma}\Gamma(\gamma+1)\right)^{-2}}{\sqrt{2\pi}\epsilon^{d}}\int_{0}^{r}g(z)dz,
}
where
\[
g(z)=z^{2\gamma+1}{}_1F_{2}\left(\gamma+\frac{1}{2}; \gamma+1, 2\gamma+1; -z^2\right)\ {}_{1}F_2\left( \lambda; \lambda+\frac{\mu}{2}, \lambda+\frac{\mu+1}{2}; -\left(\frac{z}{2\epsilon}\right)^2\right).
\]

 The product of two hypergeometric functions  ${}_{1}F_2$ can be rewritten as Kamp\'e der F\'eriet function, \citep{MGH} (1.3.30) so that 
 
\begin{equation} \begin{aligned} 
g(z)&=z^{2\gamma+1 }F \begin{smallmatrix}0  : 1;1\\0  :  2;2\\
\end{smallmatrix} \left(\begin{array}{lll}  &: \gamma+\frac12&; \lambda\\
 &: \gamma+1, 2\gamma+1&; \lambda+ \frac{\mu}{2},\lambda+\frac{\mu+1}{2}
\end{array};\-z^2,\ -\left(\frac{z}{2\epsilon}\right)^2 \right),
\end{aligned} \end{equation}
 where 
 \begin{equation}\begin{aligned}
 F \begin{smallmatrix}\mu  : \nu;\nu\\\rho  :  \sigma;\sigma\\
\end{smallmatrix} \left(\begin{array}{lll } a_1, \ldots, a_{\mu} &: b_1,\ldots,b_{\nu} &; b'_{1},\ldots,b'_{\nu} \\
c_{1},\ldots, c_{\rho} &: d_{1},\ldots,d_{\sigma}&; d'_1,\ldots, d'_{\sigma}
\end{array};\ x,\ y \right)\\
\quad =\sum_{m,n =0}^{\infty} \frac{\prod_{\ell=1}^{\mu}(a_{\ell})_{n+m}\prod_{\ell=1}^{\nu}(b_{\ell})_{n}\prod_{\ell=1}^{\nu}(b'_{\ell})_{m}}{\prod_{\ell=1}^{\rho}(c_{\ell})_{n+m}\prod_{\ell=1}^{\sigma}(d_{\ell})_{n}\prod_{\ell=1}^{\sigma}(d'_{\ell})_{m}}\frac{x^ny^m}{n!m!},
\end{aligned}
 \end{equation}
 
 where convergence is guaranteed for any  $x$ and $y$  provided that
 \[
 \mu+\nu< \rho + \sigma +1\quad {\rm{and}}\quad |x|<\infty, |y|< \infty.
 \]
 We observe that this condition is satisfied in our case where we have

 \EQ{grep}{
 g(z) = \sum_{m,n =0}^{\infty} \frac{\left(\gamma+\frac{1}{2}\right)_{n}(\lambda)_{m}}{(\gamma+1)_{n}(2\gamma+1)_{n}\left(\lambda+\frac{\mu}{2}\right)_{m}\left(\lambda+\frac{\mu+1}{2}\right)_{m}}\frac{(-1)^{m+n}}{(2\epsilon)^{2m}}\frac{z^{2(m+n+\gamma)+1}}{m!n!}.
 }
 
 Since  the pochhammer  symbol multipliers in (\ref{grep}) are all positive, the series is absolutely convergent for $|z|<\infty$ and so we can integrate term-wise to yield:
 
\[
\begin{aligned}
&\int_{0}^{r}g(z)dz =\\
&\sum_{m,n =0}^{\infty}\frac{\left(\gamma+\frac{1}{2}\right)_{n}(\lambda)_{m}}{(\gamma+1)_{n}(2\gamma+1)_{n}\left(\lambda+\frac{\mu}{2}\right)_{m}\left(\lambda+\frac{\mu+1}{2}\right)_{m}}\frac{(-1)^{m+n}}{(2\epsilon)^{2m}}\frac{r^{2(m+n+\gamma+1)}}{2(m+n+\gamma+1)m!n!}\\
&=\frac{r^{2(\gamma+1)}}{2(\gamma+1)} \sum_{m,n =0}^{\infty} \frac{(\gamma+1)_{m+n}\left(\gamma+\frac{1}{2}\right)_{n}(\lambda)_{m}(-r^{2})^n\left(-\frac{r^2}{4\epsilon^2}\right)^{m}}{(\gamma+2)_{m+n}(\gamma+1)_{n}(2\gamma+1)_{n}\left(\lambda+\frac{\mu}{2}\right)_{m}\left(\lambda+\frac{\mu+1}{2}\right)_{m}m!n!}\\
&=\frac{r^{2(\gamma+1)}}{2(\gamma+1)} F \begin{smallmatrix}1  : 1;1\\1  :  2;2
\end{smallmatrix} \left(\begin{array}{lll}  \gamma+1 &: \gamma+\frac12&; \lambda;\\
 \gamma+2 &: \gamma+1, 2\gamma+1&; \frac{\mu}{2}+\lambda,\frac{\mu}{2}+\frac{1}{2} +\lambda;
\end{array} -r^2,\ -\frac{r^{2}}{4\epsilon^{2}} \right).
\end{aligned}
\]

Consider the following identity for the Pochhamer symbol 
\EQ{pochsum}{
(a)_{m+n}=(a)_{m}(a+m)_{n}=(a)_{n}(a+m)_{n}.
}
We can use the first equality to write the above expression for the finite integral as 
\[
\begin{aligned}
&\int_{0}^{r}g(z)dz =\frac{r^{2(\gamma+1)}}{2(\gamma+1)} \\
&\times \sum_{m,n =0}^{\infty} \frac{(\gamma+1)_{m}(\gamma+1+m)_{n}\left(\gamma+\frac{1}{2}\right)_{n}(\lambda)_{m}(-r^{2})^n\left(-\frac{r^2}{4\epsilon^2}\right)^{m}}{(\gamma+2)_{m}(\gamma+2+m)_{n}(\gamma+1)_{n}(2\gamma+1)_{n}\left(\lambda+\frac{\mu}{2}\right)_{m}\left(\lambda+\frac{\mu+1}{2}\right)_{m}m!n!}\\
&=\frac{r^{2(\gamma+1)}}{2(\gamma+1)} \\
&\times \sum_{m =0}^{\infty} \frac{(\gamma+1)_{m}(\lambda)_{m}\left(-\frac{r^2}{4\epsilon^2}\right)^{m}}{(\gamma+2)_{m}\left(\lambda+\frac{\mu}{2}\right)_{m}\left(\lambda+\frac{\mu+1}{2}\right)_{m}m!}
\sum_{n=0}^{\infty}\frac{(\gamma+1+m)_{n}\left(\gamma+\frac{1}{2}\right)_{n}(-r^{2})^n}{(\gamma+2+m)_{n}(\gamma+1)_{n}(2\gamma+1)_{n}n!}.
\end{aligned}
\]

We recognise the infinite sum indexed by $n$ as a  ${}_{2}F_3$ hypergeometric series and so we can write
\EQ{first2}{
\begin{aligned}
\int_{0}^{r}g(z)dz =&\frac{r^{2(\gamma+1)}}{2(\gamma+1)} \sum_{\ell =0}^{\infty} \frac{(\gamma+1)_{\ell}(\lambda)_{\ell}\left(-\frac{r^2}{4\epsilon^2}\right)^{\ell}}{(\gamma+2)_{\ell}\left(\lambda+\frac{\mu}{2}\right)_{\ell}\left(\lambda+\frac{\mu+1}{2}\right)_{\ell}\ell!}\\
&\times{}_{2}F_3\left(\gamma+1+\ell,\gamma+\frac{1}{2};\gamma+2+\ell,\gamma+1,2\gamma+1;-r^{2}\right),
\end{aligned}
}

where, for convenience, we have switched the index from $m$ to $\ell$ in the final line. The result now follows from  equation (\ref{gammaintegral}).
\end{proof}

The closed from representation enables us to prove the main result of this section. 
\begin{proof}[Proof of Theorem \ref{Coeff}]
We recall that the spherical coefficients of the generalised Wendland functions are given by (\ref{GWcoefs}). Thus, the next stage of the investigation involves computing  (\ref{first}) 
in the limit as $r\to \infty.$ To achieve this we require an asymptotic expansion for the ${}_{2}F_3$ hypergeometric function.  An asymptotic expansion for generalised hypergeometric functions, of the kind we are interested in, is given in  \citep{Miller1997}, Equation (2.1a), it was originally derived by C.S. Meijer and a more detailed version is proven in \citep{Luke1969}.
It states that for $\vert r \vert \rightarrow \infty$, $\vert \operatorname{arg}(r)\vert \leq \frac{\pi}{2}$, 
\EQ{assgen}{
\begin{aligned}
{}_pF_{p+1}&\left( a_1,\ldots,a_p\ ;\ b_1,\ldots,b_{p+1}; -r^2\right)=\frac{\Gamma(\mathbf{b})}{\Gamma(\mathbf{a})}\sum_{k=1}^p \frac{\Gamma(a_k)\Gamma(\mathbf{a}^*-a_k)}{\Gamma(\mathbf{b}-a_k)r^{2a_k}}\\
&\times {}_{p+2}F_{p-1}\left( a_k,1+ a_k-\mathbf{b};1+a_k-\mathbf{a}^*; -\frac{1}{r^2}\right)\\
&+\frac{\Gamma(\mathbf{b})}{\Gamma(\frac12)\Gamma(\mathbf{a})}\left( \frac{1}{r^2}\right)^n\left[ 1+\bigO\left( \frac{1}{r^2}\right) \right]\cos(\xi(r)), 
\end{aligned} }

where $\mathbf{a}=(a_1,\ldots,a_p)$, $\mathbf{b}=(b_1,\ldots,b_{p+1})$,
\begin{equation*}
n=\frac12 \left( \sum_{k=1}^{p+1}b_{k}-\sum_{k=1}^{p}a_k-\frac12 \right) , \ \xi(r)=2r-\pi n+\bigO\left(\frac{1}{r}\right)
\end{equation*}
and 
\[
\begin{aligned}
\Gamma(\mathbf{a})&:=\prod_{j=1}^{p}\Gamma(a_{j}), \quad\quad  \Gamma(\mathbf{a}^*-a_k):=\prod_{j=1(j\ne k)}^{p}\Gamma(a_{j}-a_{k}), \\
\Gamma(\mathbf{b})&:=\prod_{j=1}^{p+1}\Gamma(b_{j}),\quad \quad \Gamma\left(\mathbf{b}-a_{k}\right):=\prod_{j=1}^{p+1}\Gamma(b_{j}-a_{k}).\\
 \end{aligned}
\]
We note here that  if none of the parameters appearing in the ${}_{p+2}F_{p-1}$ are negative integers, then this part can be estimated,  for $\vert r \vert \rightarrow \infty,$ by   

\begin{equation}
{}_{p+2}F_{p-1}\left( a_k, 1+a_k-\mathbf{b};1+a_k-\mathbf{a}^*; -\frac{1}{r^2}\right) =1+\bigO\left(\frac{1}{r^2}\right).
\end{equation}

An  application of (\ref{assgen})  with $p=2$ to the hypergeometric function appearing in (\ref{first}) yields an expansion with three terms
\[
{}_2F_3 \left( \
 \gamma+\ell+1,\gamma+ \frac12;\  \gamma+2+\ell,\gamma+1,2\gamma+1;
\ -r^2\right) =\sum_{j=1}^{3}T_{\ell}^{(j)}(r),
\]
where
\[
\begin{aligned}
T_{\ell}^{(1)}(r) =& \frac{\Gamma(-\ell- \frac12)}{\Gamma(\gamma+\frac12)}\frac{\Gamma(\gamma+2+\ell)\Gamma(\gamma+1)\Gamma(2\gamma+1)}{\Gamma(1)\Gamma(-\ell)\Gamma(\gamma-\ell)r^{2(\gamma+\ell+1)} }\\
& \times {}_4F_1\left(  \gamma+\ell+1, 0, \ell+1, -\gamma+\ell+1; 
\ell+\frac32;\ -\frac{1}{r^2}\right),
\end{aligned}
\]
\[
\begin{aligned}
T_{\ell}^{(2)}(r) =& \frac{\Gamma(\ell+ \frac12)}{\Gamma(\gamma+\ell+1)}\frac{\Gamma(\gamma+2+\ell)\Gamma(\gamma+1)\Gamma(2\gamma+1)}{\Gamma(\ell+ \frac32)\Gamma(\frac12)\Gamma(\gamma+ \frac12)r^{2\gamma+1}} \\
&\times {}_4F_1\left( \gamma+ \frac12, -\frac12-\ell,  \frac12,-\gamma+\frac12; \
 \frac12-\ell ;\ -\frac{1}{r^2}\right)
\end{aligned}
\]
and
\[
\begin{aligned}
T_{\ell}^{(3)}(r) =& \frac{\Gamma(\gamma+2+\ell)\Gamma(\gamma+1)\Gamma(2\gamma+1)}{\Gamma(1/2)\Gamma(\gamma+\ell+1)\Gamma(\gamma+ \frac12)r^{2\gamma+1}} \left[ 1+ O \left(\frac{1}{r^2}\right)\right] \cos(\xi(r)),\\
\xi(r)&=\left( 2r- \pi\left(\gamma+\frac{1}{2}\right)+ O\left(\frac{1}{r}\right) \right).
\end{aligned}
\]

We observe that $T_{\ell}^{(1)}(r)$ is zero because the denominator contains a Gamma function taking a non-positive integer $-\ell$ as its argument. Also,   the hypergeometric function component of $T_{\ell}^{(1)}(r)$  is equal to one this being a consequence of the zero in the list of parameters.  An application of the duplication formula and the Gamma function property $\Gamma(z+1)=z\Gamma(z)$ allows us to simplify the remaining expressions to 
\[
T_{\ell}^{(2)}(r)=\frac{2(\gamma+1+\ell)}{(2\ell+1)r^{2\gamma+1}}\frac{(2^{\gamma}\Gamma(\gamma+1))^{2}}{\pi}{}_4F_1\left( \gamma+ \frac12, -\frac12-\ell, \frac12,-\gamma+\frac12; \
 \frac12-\ell ;\ -\frac{1}{r^2}\right)
 \]
 and
\EQ{T3l}{
T_{\ell}^{(3)}(r)=\frac{(\gamma+1+\ell)}{r^{2\gamma+1}}\frac{(2^{\gamma}\Gamma(\gamma+1))^{2}}{\pi}\left[ 1+ O \left(\frac{1}{r^2}\right)\right] \cos(\xi(r)).
 }

In view of the above findings, we can return to identity (\ref{first2}) for the finite integral under investigation and write it as the sum of two components
\[
\int_{0}^{r}g(z)dz = S_{2}(r)+S_{3}(r)
\]
where
\[
S_{2}(r)=\frac{r^{2(\gamma+1)}}{2(\gamma+1)} \sum_{\ell =0}^{\infty} \frac{(\gamma+1)_{\ell}(\lambda)_{\ell}\left(-\frac{r^2}{4\epsilon^2}\right)^{\ell}}{(\gamma+2)_{\ell}\left(\lambda+\frac{\mu}{2}\right)_{\ell}\left(\lambda+\frac{\mu+1}{2}\right)_{\ell}\ell!}T_{\ell}^{(2)}(r)
\]
and
\EQ{S3def}{
S_{3}(r)=\frac{r^{2(\gamma+1)}}{2(\gamma+1)} \sum_{\ell =0}^{\infty} \frac{(\gamma+1)_{\ell}(\lambda)_{\ell}\left(-\frac{r^2}{4\epsilon^2}\right)^{\ell}}{(\gamma+2)_{\ell}\left(\lambda+\frac{\mu}{2}\right)_{\ell}\left(\lambda+\frac{\mu+1}{2}\right)_{\ell}\ell!}T_{\ell}^{(3)}(r).
}

We start with a simplification of $S_2(r)$.

\EQ{S2simp1}
{
\begin{aligned}
&S_{2}(r)=\frac{r^{2(\gamma+1)}}{2(\gamma+1)} \frac{(2^{\gamma}\Gamma(\gamma+1))^{2}}{\pi}\sum_{\ell =0}^{\infty} \frac{(\gamma+1)_{\ell}(\lambda)_{\ell}\left(-\frac{r^2}{4\epsilon^2}\right)^{\ell}}{(\gamma+2)_{\ell}\left(\lambda+\frac{\mu}{2}\right)_{\ell}\left(\lambda+\frac{\mu+1}{2}\right)_{\ell}\ell!}\frac{2(\gamma+1+\ell)}{(2\ell+1)r^{2\gamma+1}}\\
&\times {}_4F_1\left( \gamma+ \frac12, -\frac12-\ell,  \frac12,-\gamma+\frac12; \
 \frac12-\ell ;\ -\frac{1}{r^2}\right)\\
 &=\frac{(2^{\gamma}\Gamma(\gamma+1))^{2}r}{2\pi}\sum_{\ell =0}^{\infty}  \frac{(\lambda)_{\ell}\left(-\frac{r^2}{4\epsilon^2}\right)^{\ell}}{\left(\lambda+\frac{\mu}{2}\right)_{\ell}\left(\lambda+\frac{\mu+1}{2}\right)_{\ell}\ell!\left(\ell+\frac{1}{2}\right)}\\
 &\times \sum_{j=0}^{\infty}\frac{\left(\frac{1}{2}-\gamma\right)_{j}\left(\frac{1}{2}\right)_{j}\left(\frac{1}{2}+\gamma\right)_{j}\left(-\frac{1}{2}-\ell\right)_{j}}{\left(\frac{1}{2}-\ell\right)_{j}j!} \left(\frac{-1}{r^2}\right)^{j}\\
 &=\frac{(2^{\gamma}\Gamma(\gamma+1))^{2}}{2\pi}\sum_{j=0}^{\infty}\frac{(-1)^j\left(\frac{1}{2}-\gamma\right)_{j}\left(\frac{1}{2}\right)_{j}\left(\frac{1}{2}+\gamma\right)_{j}}{\left(\frac{1}{2}-j\right)j!}\frac{1}{r^{2j-1}}\\
 &\times \sum_{\ell =0}^{\infty}  \frac{\left(\frac{1}{2}-j\right)_{\ell}(\lambda)_{\ell}\left(-\frac{r^2}{4\epsilon^2}\right)^{\ell}}{\left(\lambda+\frac{\mu}{2}\right)_{\ell}\left(\lambda+\frac{\mu+1}{2}\right)_{\ell}\left(\frac{3}{2}-j\right)_{\ell}\ell!}\\
 &=\frac{(2^{\gamma}\Gamma(\gamma+1))^{2}}{2\pi}\sum_{j=0}^{\infty}\frac{(-1)^{j}\left(\frac{1}{2}-\gamma\right)_{j}\left(\frac{1}{2}\right)_{j}\left(\frac{1}{2}+\gamma\right)_{j}}{\left(\frac{1}{2}-j\right)j!}\\
 &\times \frac{1}{r^{2j-1}}  {}_{2}F_{3}\left(\frac{1}{2}-j,\lambda;\lambda+\frac{\mu}{2},\lambda+\frac{\mu+1}{2},\frac{3}{2}-j;-\frac{r^2}{4\epsilon^2}\right),
\end{aligned}
}

where, in the first and the second equations we have used, respectively, the following identities  
\EQ{ids}{
\frac{(\gamma+1)_{\ell}}{(\gamma+2)_{\ell}}=\frac{\gamma+1}{\gamma+1+\ell}\quad {\rm{and}}\quad \frac{1}{\ell+\frac{1}{2}}\frac{\left(-\frac{1}{2}-\ell\right)_{j}}{\left(\frac{1}{2}-\ell\right)_{j}}=
\frac{1}{\frac{1}{2}-j}\frac{\left(\frac{1}{2}-j\right)_{\ell}}{\left(\frac{3}{2}-j\right)_{\ell}}.
}
In the final line we recognise the series indexed by $\ell$ as a $ {}_{2}F_{3}$ hypergeometric function. Applying the asymptotic expansion, as in the earlier development,  we find that
\[
 _{2}F_{3}\left(\frac{1}{2}-j,\lambda;\lambda+\frac{\mu}{2},\lambda+\frac{\mu+1}{2},\frac{3}{2}-j;-\frac{r^2}{4\epsilon^2}\right)=\sum_{k=1}^{3}\tilde{T}_{j}^{(k)}(r).
 \]
 The first contribution from the sum is
  \[
 \begin{aligned}
 \tilde{T}_{j}^{(1)}(r)&=\frac{\Gamma\left(\lambda+\frac{\mu}{2}\right)\Gamma\left(\lambda+\frac{\mu+1}{2}\right)\Gamma\left(\frac{3}{2}-j\right)\Gamma\left(\lambda+j-\frac{1}{2}\right)}{\Gamma\left(\lambda\right)\Gamma\left(\lambda+\frac{\mu}{2}+j-\frac{1}{2}\right)\Gamma\left(\lambda+\frac{\mu}{2}+j\right)}\left(\frac{2\epsilon}{r}\right)^{2\left(\frac{1}{2}-j\right)}\\
 &\times {}_4F_1\left(\frac{1}{2}-j,0,\frac{3-\mu}{2} -j-\lambda,\frac{2-\mu}{2} -j-\lambda; -\frac{4\epsilon^2}{r^2}\right)\\
 &=\frac{\Gamma\left(\lambda-\frac{1}{2}\right)\Gamma\left(\lambda+\frac{\mu+1}{2}\right)\Gamma\left(\frac{3}{2}-j\right)\left(\lambda-\frac{1}{2}\right)_{j}}{\Gamma\left(\lambda\right)\Gamma\left(\lambda+\frac{\mu-1}{2}\right)\left(\lambda+\frac{\mu-1}{2}\right)_{j}\left(\lambda+\frac{\mu}{2}\right)_{j}}\left(\frac{r}{2\epsilon}\right)^{2j-1},
 \end{aligned}
 \]
 where, in the final line we have used the observation that the ${}_{4}F_1$ hypergeometric function is equal to $1$ due to the zero appearing as one of its coefficients. Additionally, we have used Pochammer symbol notation as this will become useful later in the development.
  We observe that
 \EQ{obs1}{
 \begin{aligned}
& \lim_{r\to \infty}\frac{\tilde{T}_{j}^{(1)}(r)}{r^{2j-1}}=
\frac{\Gamma\left(\lambda-\frac{1}{2}\right)\Gamma\left(\lambda+\frac{\mu+1}{2}\right)\Gamma\left(\frac{3}{2}-j\right)\left(\lambda-\frac{1}{2}\right)_{j}}{\Gamma\left(\lambda\right)\Gamma\left(\lambda+\frac{\mu-1}{2}\right)\left(\lambda+\frac{\mu-1}{2}\right)_{j}\left(\lambda+\frac{\mu}{2}\right)_{j}\left(2\epsilon\right)^{2j-1}}\\
&=\frac{\Gamma\left(\lambda-\frac{1}{2}\right)\left(2\lambda+\mu-1\right)\epsilon}{\Gamma\left(\lambda\right)}\frac{\Gamma\left(\frac{3}{2}-j\right)\left(\lambda-\frac{1}{2}\right)_{j}}{\left(\lambda+\frac{\mu-1}{2}\right)_{j}\left(\lambda+\frac{\mu}{2}\right)_{j}\left(4\epsilon^2\right)^{j}}\,\,\, \,\,\quad\quad {\rm{for}}\,\, j\ge0.
\end{aligned}
 }
 
  The second contribution from the sum is
 \[
 \begin{aligned}
 \tilde{T}_{j}^{(2)}(r)&=\frac{\Gamma\left(\lambda+\frac{\mu}{2}\right)\Gamma\left(\lambda+\frac{\mu+1}{2}\right)\left(\frac{3}{2}-j\right)}{\Gamma\left(\frac{\mu}{2}\right)\Gamma\left(\frac{\mu+1}{2}\right)\left(\frac{3}{2}-\lambda-j\right)}\left(\frac{2\epsilon}{r}\right)^{2\lambda}\\
 &\times {}_4F_1\left(\lambda,1-\frac{\mu}{2}, 1-\frac{\mu+1}{2},\lambda+j-\frac{1}{2};\lambda+j+\frac{1}{2}; -\frac{4\epsilon^2}{r^2}\right).\\
 \end{aligned}
 \]
 We note that if either $\mu=1$ or $\mu=2$ then $_4F_1$ hypergeometric function appearing in the formula above collapses to $1$ again, due to a zero appearing as one of its coefficients in these cases. More generally, one can deduce that
 \[
 {}_4F_1\left(\lambda,1-\frac{\mu}{2}, 1-\frac{\mu+1}{2},\lambda+j-\frac{1}{2};\lambda+j+\frac{1}{2}; -\frac{4\epsilon^2}{r^2}\right)= 1+O\left(r^{-2}\right)
  \]
  and so conclude, for the second contribution, that
\[   \tilde{T}_{j}^{(2)}(r)=\frac{\Gamma\left(\lambda+\frac{\mu}{2}\right)\Gamma\left(\lambda+\frac{\mu+1}{2}\right)\left(\frac{3}{2}-j\right)}{\Gamma\left(\frac{\mu}{2}\right)\Gamma\left(\frac{\mu+1}{2}\right)\left(\frac{3}{2}-\lambda-j\right)}\left(\frac{2\epsilon}{r}\right)^{2\lambda}\Bigl(1+O\left(r^{-2}\right)\Bigr).
 \] 
 
 Given that $\lambda >\frac{1}{2}$ we observe that
 \EQ{obs2}{ \lim_{r\to \infty}\frac{1}{r^{2j-1}}\tilde{T}_{j}^{(2)}(r)=0 \quad {\rm{for}}\,\, j\ge0.
 }

The third contribution from the sum is
\[
 \tilde{T}_{j}^{(3)}(r)=\frac{\Gamma\left(\lambda+\frac{\mu}{2}\right)\Gamma\left(\lambda+\frac{\mu+1}{2}\right)\Gamma\left(\frac{3}{2}-j\right)}
 {\sqrt{\pi}\Gamma\left(\frac{1}{2}-j\right)\Gamma\left(\lambda\right)}\frac{(2\epsilon)^{\lambda+\mu+1}}{r^{\lambda+\mu+1}}\Bigl(1+O\left(r^{-2}\right)\Bigr)\cos\left(\xi(r)\right),
 \]
 and, in this case,  since $\lambda, \mu >0$ we observe that
  \EQ{obs3}{ \lim_{r\to \infty}\frac{1}{r^{2j-1}}\tilde{T}_{j}^{(3)}(r)=0 \quad {\rm{for}}\,\, j\ge0.
 }
 
 Next we turn to the simplification of $S_{3}(r)$ given in (\ref{S3def}) for which, using the first identity of (\ref{ids}) and the estimate (\ref{T3l}), we have
\[
\begin{aligned}
S_{3}(r)&=
\frac{\left(2^{\gamma}\Gamma\left(\gamma+1\right) \right)^{2}r}{\pi}\sum_{\ell=0}^{\infty}\frac{(\lambda)_{\ell}}{\left(\lambda+\frac{\mu}{2}\right)_{\ell}\left(\lambda+\frac{\mu+1}{2}\right)_{\ell}\ell!}\left(-\frac{r^2}{4\epsilon^2}\right)^{\ell}\Bigl( 1+ O \left(r^{-2}\right)\Bigr) \cos(\xi(r))
\\
&=\frac{\left(2^{\gamma}\Gamma\left(\gamma+1\right) \right)^{2}r}{\pi}{}_1F_2\left( \lambda; \lambda+\frac{\mu}{2},\lambda+\frac{\mu+1}{2}
 ;\ -\frac{r^2}{4\epsilon^2}\right)\Bigl( 1+ O \left(r^{-2}\right)\Bigr) \cos(\xi(r)).
\end{aligned}
\]

We have encountered this ${}_1F_2$ hypergeometric function earlier in the paper, specifically (\ref{ass1F2}) where one can see that it decays to zero, as $r\to \infty,$ at a rate that is faster than $\frac{1}{r}.$ Thus we can conclude that

 \EQ{obs4}
 {
 \lim_{r\to \infty}S_{3}(r)=0.
 }
 
 Now,  using the four observations (\ref{obs1})-(\ref{obs4}) we can conclude that
 
 
 
\EQ{obs5}
{
\begin{aligned}
  \lim_{r\to \infty}\int_{0}^{r}g(z)dz&=\frac{\left(2^{\gamma}\Gamma\left(\gamma+1\right) \right)^{2}}{2\pi} \frac{\Gamma\left(\lambda-\frac{1}{2}\right)\left(2\lambda+\mu-1\right)\epsilon}{\Gamma\left(\lambda\right)} \\
&\times \sum_{j=0}^{\infty}\frac{(-1)^{j}\left(\frac{1}{2}-\gamma\right)_{j}\left(\frac{1}{2}+\gamma\right)_{j}\left(\frac{1}{2}\right)_{j}\Gamma\left(\frac{3}{2}-j\right)\left(\lambda-\frac{1}{2}\right)_{j}}{\left(\frac{1}{2}-j\right)j!\left(\lambda+\frac{\mu-1}{2}\right)_{j}\left(\lambda+\frac{\mu}{2}\right)_{j}(4\epsilon^2)^{j}}\\
&=\frac{\left(2^{\gamma}\Gamma\left(\gamma+1\right) \right)^{2}}{2\sqrt{\pi}}\frac{\Gamma\left(\lambda-\frac{1}{2}\right)\left(2\lambda+\mu-1\right)\epsilon}{\Gamma\left(\lambda\right)} \\
&\times \sum_{j=0}^{\infty}\frac{\left(\frac{1}{2}-\gamma\right)_{j}\left(\frac{1}{2}+\gamma\right)_{j}\left(\lambda-\frac{1}{2}\right)_{j}}{j!\left(\lambda+\frac{\mu-1}{2}\right)_{j}\left(\lambda+\frac{\mu}{2}\right)_{j}(4\epsilon^2)^{j}},
\end{aligned}
}
where we have used $\Gamma\left(\frac{3}{2}-j\right)=\left(\frac{1}{2}-j\right)\Gamma\left(\frac{1}{2}-j\right)$ and the identity $\left(\frac{1}{2}\right)_{j}\Gamma\left(\frac{1}{2}-j\right)=\sqrt{\pi}(-1)^{j}.$ The infinite series above can be represented as a ${}_3F_2$ hyperpgeometric function and  so we can conclude that 
\[
\begin{aligned}
\lim_{r\to \infty}\int_{0}^{r}g(z)dz=&\frac{\left(2^{\gamma}\Gamma\left(\gamma+1\right) \right)^{2} \Gamma\left(\lambda-\frac{1}{2}\right)\left(2\lambda+\mu-1\right)\epsilon}
{2\sqrt{\pi}\Gamma\left(\lambda\right)}\\
&\times
{}_3F_2\left(\frac{1}{2}-\gamma,\frac{1}{2}+\gamma,\lambda-\frac{1}{2};\lambda+\frac{\mu-1}{2},\lambda+\frac{\mu}{2};
\frac{1}{4\epsilon^2}\right).
\end{aligned}
\]
 Revisiting the integral (\ref{gammaintegral}) we can use the above to deduce that $\lim_{r\to \infty}I_{\gamma}(r)=$
\[
\begin{aligned}
& \frac{C_{\lambda,\mu}\Gamma\left(\lambda-\frac{1}{2}\right)(2\lambda+\mu-1)}{2\sqrt{\pi}\epsilon^{d-1}\Gamma\left(\lambda\right)}
{}_3F_2\left(\frac{1}{2}-\gamma,\frac{1}{2}+\gamma,\lambda-\frac{1}{2};\lambda+\frac{\mu-1}{2},\lambda+\frac{\mu}{2};
\frac{1}{4\epsilon^2}\right)\\
&= (2\pi)^{\frac{d-1}{2}}\frac{C_{\lambda-\frac{1}{2},\mu}}{\sqrt{2\pi}\epsilon^{d-1}}{}_3F_2\left(\frac{1}{2}-\gamma,\frac{1}{2}+\gamma,\lambda-\frac{1}{2};\lambda+\frac{\mu-1}{2},\lambda+\frac{\mu}{2};
\frac{1}{4\epsilon^2}\right),
\end{aligned}
\]
where the final equality can be deduced from the definition of the constant $C_{\lambda,\mu}$ (\ref{ClamMu}).
In particular, setting $\gamma = m+\frac{d-2}{2},$ we complete the proof.  
\end{proof}
\section{Asymptotic decay of  $\widehat{\psi_{\mu,\alpha}^{(\epsilon)}}(m)$}
We can investigate the asymptotic behaviour of the spherical Fourier coefficients of the generalised Wendland function by focusing upon the hypergeometric series appearing in (\ref{swc}).

\theo{Theorem}{Asdecay}{
The spherical Fourier coefficients of the generalised Wendland functions $\Psi_{\mu,\alpha}^{(\epsilon)}(\boldsymbol\xi,\boldsymbol\eta)$ defined in  \eqref{reskernel} exhibit the following precise asymptotic decay
\EQ{decaywendsphere}
{
\widehat{\psi_{\mu,\alpha}^{(\epsilon)}}(m)\sim (2\pi)^{\frac{d-1}{2}}\frac{2^{\lambda-\frac{1}{2}}\Gamma\left(\lambda-\frac{1}{2}\right)\mu}{\sqrt{2\pi}}\frac{\epsilon^{2\alpha+1}}{\left(m+\frac{d-1}{2}\right)^{2\lambda-1}}.
}
As in the Euclidean case, this asymptotic behaviour implies that there 
 exist positive constants $c_{1}<c_{2}$ such that (\ref{decaysp}) holds and thus the zonal kernel $\Psi_{\mu,\alpha}^{(\epsilon)}(\boldsymbol\xi,\boldsymbol\eta)$  is reproducing for the Sobolev space $H^{\lambda-\frac{1}{2}}(S^{d-1}).$

}

\begin{proof}
 We begin by separating out two cases. First, when $d\ge 2$ is odd then, due to the appearance of a negative integer in the coefficient list, the series terminates and so it collapses to a hypergeometric polynomial of the form
\EQ{genhypp}{
{_3F_2}\left(-n,n+c, a\,\,;\,\,b_{1}, b_{2}\,\,;\,\, z\right),\,\,\,\,\, n \in \Z_{+}.
}
In the second case, where $d \ge 2$ is even, the series does not terminate.

A survey of the literature in this area shows that there are very few known asymptotic results that apply to general ${_3F_2}$ hypergeometric functions for large parameters. Most of the known results apply to the case where the series terminates and, in this regard,
 we are fortunate that the limiting behaviour of  (\ref{genhypp})
 as $n\to \infty$  is covered in \citep{F3} where it is shown that provided none of the hypergeometric parameters $(a, b_{1}, b_{2}, c)$ in (\ref{genhypp}) coincide with zero or with a negative integer, then  
with the following definitions
\EQ{alpha}
{2\alpha = a-b_{1}-b_{2}+\frac{1}{2}\,\,\,\,\,\,\,{\rm{and}}\,\,\,\,\,\,\,\, z = \sin^{2}(\theta/2)\in (0,1),}
we have the following asymptotic results

$\diamond$ For $z \in (0,1):$

\EQ{crucialz}{
\begin{aligned}
&{_3F_2}\left(-n,n+c, a\,;\,b_{1},b_{2}\,;\, z\right)=\frac{\Gamma(b_{1})\Gamma(b_{2})}{\Gamma(b_{1}-a)\Gamma(b_{2}-a)}\frac{1}{z^{a}}\frac{1}{(n+c)^{2a}}\left[1+O\left(\frac{1}{(n+c)}\right)\right]\\
&\\
\, & \\
&+\frac{\Gamma(b_{1})\Gamma(b_{2})(n+c)^{2\alpha}}{\sqrt{\pi}\Gamma(a)}\left[\frac{\left( \sin^{2}(\theta/2)\right)^{\alpha}\cos\left(\left(n+\frac{c}{2}\right)\theta +\pi \alpha\right) }{\left( \cos^{2}(\theta/2)\right)^{\alpha+\frac{c}{2}}}\right]
+O\left((n+c)^{2\alpha-1}\right).
\end{aligned}
}

$\diamond$ For $z=1:$


\EQ{crucialone}{
\begin{aligned}
&{_3F_2}\left(-n,n+c, a\,;\,b_{1},b_{2}\,;\, 1\right)=\frac{\Gamma(b_{1})\Gamma(b_{2})}{\Gamma(b_{1}-a)\Gamma(b_{2}-a)}\frac{1}{(n+c)^{2a}}\left[1+O\left(\frac{1}{(n+c)}\right)\right]+\\
\, & \\
&\frac{(-1)^{n}\Gamma(b_{1})\Gamma(b_{2})\Gamma(n+2c+4\alpha)}{\Gamma\left(c+2\alpha+\frac{1}{2}\right)\Gamma(a)\Gamma (n+c)}\left[1-\frac{(c+4\alpha)\left(c+2\alpha-\frac{1}{2}\right)}{(n+c)}+O\left(\frac{1}{(n+c)^{2}}\right)\right].
\end{aligned}
}

In our case, $n= m+\frac{d-3}{2}, c = 1, a = \lambda-\frac{1}{2}, b_{1}=\lambda+\frac{\mu-1}{2},$ $b_{2}=\lambda+\frac{\mu}{2}$ and $z = \frac{1}{4\epsilon^{2}}.$ Since, $\lambda,\mu>\frac{1}{2}$ it is straight-forward to check that these parameters satisfy the conditions associated with the above asymptotic formulae. Thus, we can set
\[
2\alpha:=-\left(\lambda+\mu-\frac{1}{2}\right)\,\,\,\, {\rm{and}}\,\,\,\, \sin^{2}\left(\frac{\theta}{2}\right) := \frac{1}{4\epsilon^2},
\]
and employ (\ref{crucialz}) and (\ref{crucialone}) to yield the following asymptotic results.

{\bf{$\diamond$}} For $\epsilon > \frac{1}{2}:$

\EQ{epsgthalf}{
\begin{aligned}
&{}_3F_2\left(- \left( m+\frac{d-3}{2}\right),m+\frac{d-1}{2},\lambda-\frac{1}{2};\lambda+\frac{\mu-1}{2},\lambda+\frac{\mu}{2};
\frac{1}{4\epsilon^2}\right)\\ \, &
\\&
=\frac{\Gamma\left(\lambda+\frac{\mu-1}{2}\right)\Gamma\left(\lambda+\frac{\mu}{2}\right)\left(2\epsilon\right)^{2\lambda-1}}{\Gamma\left(\frac{\mu}{2}\right)\Gamma\left(\frac{\mu+1}{2}\right)\left(m+\frac{d-1}{2}\right)^{2\lambda-1}}
\left[1+O\left(\frac{1}{m+\frac{d-1}{2}}\right)\right]\\
\, & \\
&+\frac{\Gamma\left(\lambda+\frac{\mu-1}{2}\right)\Gamma\left(\lambda+\frac{\mu}{2}\right)2\epsilon\left(4\epsilon^2-1\right)^{\frac{1}{2}(\lambda+\mu-\frac{3}{2})}}{\sqrt{\pi}\Gamma\left(\lambda-\frac{1}{2}\right)\left(m+\frac{d-1}{2}\right)^{\lambda+\mu-\frac{1}{2}}}
\\ \, & \\
 & \times
\left[\cos\left(\left(m+\frac{d-2}{2}\right)\theta -\frac{\pi}{2}\left(\lambda+\mu-\frac{1}{2}\right)\right)\right] +O\left(\frac{1}{\left(m+\frac{d-1}{2}\right)^{\lambda+\mu+\frac{1}{2}}}\right).
\end{aligned}
}

We observe that since $\mu \ge \lambda$ we have that 
\[
\lambda+\mu+\frac{1}{2}>\lambda+\mu-\frac{1}{2}\ge 2\lambda -\frac{1}{2}>2\lambda-1,
\]
and so, in this case, the precise asymptotic decay of the hypergeometric function is determined by the first term of (\ref{epsgthalf}) and is given by

\EQ{rep}{
\begin{aligned}
&{}_3F_2\left(- \left( m+\frac{d-3}{2}\right),m+\frac{d-1}{2},\lambda-\frac{1}{2};\lambda+\frac{\mu-1}{2},\lambda+\frac{\mu}{2};
\frac{1}{4\epsilon^2}\right)\\ \, &
\\&
\sim
\frac{\Gamma\left(\lambda+\frac{\mu-1}{2}\right)\Gamma\left(\lambda+\frac{\mu}{2}\right)\left(2\epsilon\right)^{2\lambda-1}}{\Gamma\left(\frac{\mu}{2}\right)\Gamma\left(\frac{\mu+1}{2}\right)\left(m+\frac{d-1}{2}\right)^{2\lambda-1}}
= \frac{\Gamma\left(2\lambda+\mu-1\right)\epsilon^{2\lambda-1}}{\Gamma\left(\mu\right)}\frac{1}{\left(m+\frac{d-1}{2}\right)^{2\lambda-1}}\\
&= \frac{2^{\lambda-\frac{1}{2}}\Gamma\left(\lambda-\frac{1}{2}\right)\mu}{C_{\lambda-\frac{1}{2},\mu}}\frac{\epsilon^{2\lambda-1}}{\left(m+\frac{d-1}{2}\right)^{2\lambda-1}},
\end{aligned}
}
  where the first equality follows by applying the duplication formula (\ref{gammadup}) with $z=\lambda+\frac{\mu-1}{2}$ and $z=\frac{\mu}{2},$ the final equality follows from the definition of $C_{\lambda-\frac{1}{2},\mu}$ (\ref{ClamsphMu}). The asymptotic result for the spherical coefficients  follows, in this case, directly from their definition (\ref{swc}).


$\diamond$ For $\epsilon =\frac{1}{2}:$

\EQ{epshalf}{
\begin{aligned}
&{}_3F_2\left(- \left( m+\frac{d-3}{2}\right),m+\frac{d-1}{2},\lambda-\frac{1}{2};\lambda+\frac{\mu-1}{2},\lambda+\frac{\mu}{2};
1\right)\\ \, &
\\&
=\frac{\Gamma\left(\lambda+\frac{\mu-1}{2}\right)\Gamma\left(\lambda+\frac{\mu}{2}\right)}{\Gamma\left(\frac{\mu}{2}\right)\Gamma\left(\frac{\mu+1}{2}\right)\left(m+\frac{d-1}{2}\right)^{2\lambda-1}}
\left[1+O\left(\frac{1}{m+\frac{d-1}{2}}\right)\right]\\
\, & \\
&+\frac{(-1)^{m+\frac{d-3}{2}}\Gamma\left(\lambda+\frac{\mu-1}{2}\right)\Gamma\left(\lambda+\frac{\mu}{2}\right)\Gamma\left(m+\frac{d-1}{2}-2(\lambda+\mu-1)\right)}{\Gamma\left(2-\lambda-\mu\right)\Gamma\left(\lambda-\frac{1}{2}\right)\Gamma\left(m+\frac{d-1}{2}\right)}\\
\, & \\
&\times \left[1-\frac{2(1-\lambda-\mu)^{2}}{m+\frac{d-1}{2}}+O\left(\frac{1}{\left(m+\frac{d-1}{2}\right)^{2}}\right)\right].
\end{aligned}
}

Applying the asymptotic formula for the Gamma function \citep{Abramowitz1964} equation 6.1.39,
\EQ{asgam}{
\Gamma\left(az+b\right)\sim \sqrt{2\pi}\exp(-az)(az)^{az+b-\frac{1}{2}},\quad a>0,
}
with $z= m+\frac{d-1}{2}$ we find that
\[
\frac{\Gamma\left(m+\frac{d-1}{2}-2(\lambda+\mu-1)\right)}{\Gamma\left(m+\frac{d-1}{2}\right)}\sim
\frac{1}{\left(m+\frac{d-1}{2}\right)^{2(\lambda+\mu-1)}}.
\]

Clearly, since $\mu\ge \lambda$ we have that $2(\lambda+\mu-1)\ge 2(2\lambda-1)>2\lambda-1,$ and so the asymptotic behaviour of the hypergeometric function is determined by the first term of (\ref{epshalf}). Mirroring the concluding development in (\ref{rep}) we have that
\[
\begin{aligned}
&{}_3F_2\left(- \left( m+\frac{d-3}{2}\right),m+\frac{d-1}{2},\lambda-\frac{1}{2};\lambda+\frac{\mu-1}{2},\lambda+\frac{\mu}{2};
1\right)\\ \, &
\\&
\sim\frac{\Gamma\left(\lambda+\frac{\mu-1}{2}\right)\Gamma\left(\lambda+\frac{\mu}{2}\right)}{\Gamma\left(\frac{\mu}{2}\right)\Gamma\left(\frac{\mu+1}{2}\right)\left(m+\frac{d-1}{2}\right)^{2\lambda-1}}=\frac{2^{\lambda-\frac{1}{2}}\Gamma\left(\lambda-\frac{1}{2}\right)\mu}{C_{\lambda-\frac{1}{2},\mu}}\frac{\left(\frac{1}{2}\right)^{2\lambda-1}}{\left(m+\frac{d-1}{2}\right)^{2\lambda-1}},
\end{aligned}
\]
and the asymptotic result for the spherical coefficients  follows from (\ref{swc}).

In the case where $d\ge 2$ is even we can access an asymptotic formula more directly. We have that
\[
\begin{aligned}
&{}_3F_2\left(- \left( m+\frac{d-3}{2}\right),m+\frac{d-1}{2},\lambda-\frac{1}{2};\lambda+\frac{\mu-1}{2},\lambda+\frac{\mu}{2};
\frac{1}{4\epsilon^2}\right)\\
&= \sum_{j=0}^{\infty}\frac{\left(- \left( m+\frac{d-3}{2}\right)\right)_{j}\left(m+\frac{d-1}{2}\right)_{j}\left(\lambda-\frac{1}{2}\right)_{j}}{\left(\lambda+\frac{\mu-1}{2}\right)_{j}\left(\lambda+\frac{\mu}{2}\right)_{j}j!}\left(\frac{1}{4\epsilon^2}\right)^j\\
&= \sum_{j=0}^{\infty}\frac{\Gamma\left(m+\frac{d-1}{2}+j\right)}{\Gamma\left(m+\frac{d-1}{2}-j\right)}\frac{\left(\lambda-\frac{1}{2}\right)_{j}}{\left(\lambda+\frac{\mu-1}{2}\right)_{j}\left(\lambda+\frac{\mu}{2}\right)_{j}j!}\left(\frac{-1}{4\epsilon^2}\right)^j,
\end{aligned}
\]
where, in the final line, we have used the identity
\[
(-a)_{j}=(-1)^{j}\frac{\Gamma\left(a+1\right)}{\Gamma\left(a+1-j\right)}.
\]

An application of (\ref{asgam}) yields
\[
\begin{aligned}
&{}_3F_2\left(- \left( m+\frac{d-3}{2}\right),m+\frac{d-1}{2},\lambda-\frac{1}{2};\lambda+\frac{\mu-1}{2},\lambda+\frac{\mu}{2};
\frac{1}{4\epsilon^2}\right)\\
&\sim \sum_{j=0}^{\infty}
\frac{\left(\lambda-\frac{1}{2}\right)_{j}}
{\left(\lambda+\frac{\mu-1}{2}\right)_{j}
\left(\lambda+\frac{\mu}{2}\right)_{j}j!}
\left(-\left(\frac{m+\frac{d-1}{2}}{2\epsilon}\right)^{2}\right)^j\\
&={}_1F_2\left(\lambda-\frac{1}{2};\lambda+\frac{\mu-1}{2},\lambda+\frac{\mu}{2};-\left(\frac{m+\frac{d-1}{2}}{2\epsilon}\right)^{2}\right)\\
&\sim \frac{\Gamma\left(\lambda+\frac{\mu-1}{2}\right)\Gamma\left(\lambda+\frac{\mu}{2}\right)\left(2\epsilon\right)^{2\lambda-1}}{\Gamma\left(\frac{\mu}{2}\right)\Gamma\left(\frac{\mu+1}{2}\right)\left(m+\frac{d-1}{2}\right)^{2\lambda-1}}
=  \frac{2^{\lambda-\frac{1}{2}}\Gamma\left(\lambda-\frac{1}{2}\right)\mu}{C_{\lambda-\frac{1}{2},\mu}}\frac{\epsilon^{2\lambda-1}}{\left(m+\frac{d-1}{2}\right)^{2\lambda-1}},
\end{aligned}
\]
where the asymptotic result final line follows from the asymptotic formula  (\ref{ass1F2}) for the ${}_1F_2$ hypergeometric function and the final equality is the same as in the final line of (\ref{rep}). Once again, the asymptotic result for the spherical Fourier coefficients follows from their definition (\ref{swc}).
\end{proof} 
 
We conclude the paper by drawing the reader's attention to the  close connection between the asymptotic formula for the decay of the Fourier transform of the generalised Wendland functions and that of the associated spherical Fourier coefficients. Recalling that $\lambda = \frac{d+1}{2}+\alpha$ we can define
 \[
\mathcal{K}_{\mu,\alpha}^{(d)}=\frac{2^{\frac{d+1}{2}+\alpha}\Gamma\left(\frac{d+1}{2}+\alpha\right)\mu\epsilon^{2\alpha+1}}{\sqrt{2\pi}},
 \]
 then revisiting (\ref{assdecayR}) and (\ref{decaywendsphere}) we have that
 
\[  \widehat{\phi_{\mu,\alpha}^{(\epsilon)}}(z)\sim  \frac{\mathcal{K}_{\mu,\alpha}^{(d)}}{\|z\|^{d+1+2\alpha}}\quad{\rm{and}}\quad\widehat{\psi_{\mu,\alpha}^{(\epsilon)}}(m)\sim (2\pi)^{\frac{d-1}{2}} \frac{\mathcal{K}_{\mu,\alpha}^{(d-1)}}{\left(m+\frac{d-1}{2}\right)^{d+2\alpha}}.
\]

\subsection*{Funding:}
The work of Janin Jäger was funded by the Deutsche Forschungsgemeinschaft (DFG - German research foundation) - Projektnummer: 461449252 and by the Justus-Liebig University as part of the Just'us-fellowship.

\par\bigskip\bigskip\noindent
\bibliographystyle{abbrvnat}\

\end{document}